\documentclass[12pt,a4paper]{article}
\usepackage[utf8]{inputenc}
\usepackage[english]{babel}
\usepackage[OT1]{fontenc}
\usepackage{amsmath}
\usepackage{amsfonts}
\usepackage{amssymb}
\usepackage[dvipsnames]{xcolor}
\usepackage{hyperref}
\usepackage{cleveref}
\usepackage{tikz}
\usepackage{pgfplots}
\usepackage{float}
\floatstyle{boxed}
\restylefloat{figure}

\def\1{\mathbf 1}

\def\1{\bold 1}

\usepackage{amsthm}
\theoremstyle{theorem}
\newtheorem{theorem}{Theorem}

\newtheorem{remark}[theorem]{Remark}

\theoremstyle{plain}
\newtoks\thehProclaim
\newtheorem*{Proclaim}{\the\thehProclaim}

\theoremstyle{definition}
\newtoks\thehDefinition
\newtheorem*{Definition}{\the\thehDefinition}

\providecommand{\keywords}[1]{\textbf{{Key words:}} #1}
\providecommand{\subjclass}[1]{\textit{{2010 Mathematics Subject Classification.}} #1}
\usepackage[paperwidth=20.5cm,paperheight=29.7cm,width=16.5cm,height=25.5cm,top=24mm,headsep=11mm]{geometry}

\begin{document}
\title{Note on quantitative homogenization results for~parabolic systems in $\mathbb{R}^d$
\footnote{\subjclass{Primary 35B27. Secondary 	35K45.}\break Research was supported by «Native towns», a social investment program of PJSC «Gazprom Neft», by the postdoctoral fellowship of the Hausdorff Research Institute for Mathematics in Bonn, Germany, and, last but not least, by the SPbU travel money informally called ,,100 RUB from Eremeev''.}
}
\author{Yulia Meshkova \footnote{Chebyshev Laboratory, St. Petersburg State University, 14th Line V.O., 29B, Saint Petersburg 199178 Russia.  E-mail: {\tt{y.meshkova@spbu.ru}, \tt{juliavmeshke@yandex.ru}}. }}

\maketitle

\begin{abstract}
In $L_2(\mathbb{R}^d;\mathbb{C}^n)$, we consider a semigroup $e^{-tA_\varepsilon}$, $t\geqslant 0$, generated by a matrix elliptic second order differential operator $A_\varepsilon \geqslant 0$. Coefficients of $A_\varepsilon$ are periodic, depend on $\mathbf{x}/\varepsilon$ and oscillate rapidly as $\varepsilon \rightarrow 0$. Approximations for $e^{-tA_\varepsilon}$ were obtained by T.~A.~Suslina (2004, 2010) via the spectral method and by V.~V.~Zhikov and S.~E.~Pastukhova (2006) via the shift method. In the present note, we give another short proof based on the contour integral representation for the semigroup and approximations for the resolvent with two-parametric error estimates obtained by T.~A.~Suslina (2015).
\end{abstract}


\section*{Introduction}

The subject of this note is quantitative estimates in periodic homogenization, i. e., approximations for the corresponding resolving operator in the uniform operator topology.  There are several approaches to obtaining results of such type, see \cite{BSu,CDGr,Sh,ZhPasUMN}. 

In introduction, let us consider the simplest elliptic operator $A_\varepsilon =-\mathrm{div}g(\varepsilon^{-1}\mathbf{x})\nabla$, $\varepsilon >0$, acting in $L_2(\mathbb{R}^d)$. Here $g$ is a~periodic positive definite matrix-valued function such that $g,g^{-1}\in L_\infty$. Let $u_\varepsilon$ be the solution of the equation $A_\varepsilon u_\varepsilon +u_\varepsilon =F$, where $F\in L_2(\mathbb{R}^d)$. The homogenization problem is to describe the behavior of the solution $u_\varepsilon$ in the small period limit $\varepsilon\rightarrow 0$. The classical result is that $u_\varepsilon \rightarrow u_0$ in the $L_2$-norm, where the limit function $u_0$ is the solution of the equation of the same type $A^0u_0+u_0=F$, where $A^0=-\mathrm{div}\,g^0\nabla$ is the so-called \textit{effective} operator  with the constant matrix $g^0$. 

By using the spectral method, M.~Sh.~Birman and T.~A.~Suslina \cite{BSu} proved that $\Vert u_\varepsilon -u_0\Vert _{L_2}\leqslant C\varepsilon\Vert F\Vert _{L_2}$. This estimate can be rewritten as approximation for the resolvent $(A_\varepsilon +I)^{-1}$ in the uniform operator topology. Approximations for the semigroup $e^{-tA_\varepsilon}$, $t\geqslant 0$, were obtained in \cite{Su04,Su07}, and \cite{ZhPAs_parabol} via the spectral and shift methods, respectively:
\begin{equation}
\label{Intr exp}
\Vert e^{-tA_\varepsilon }-e^{-tA^0}\Vert _{L_2(\mathbb{R}^d)\rightarrow L_2(\mathbb{R}^d)}\leqslant C\varepsilon (t+\varepsilon ^2)^{-1/2}.
\end{equation} 
Since the point zero is the lower edge of the spectra for $A_\varepsilon$ and $A^0$,  estimate \eqref{Intr exp} can be treated as a stabilization result for $t\rightarrow \infty$. 
Later in \cite{MSu}, it was observed that quantitative results for parabolic problems can be derived from  corresponding elliptic results with the help of  identity $e^{-tA_\varepsilon}=-\frac{1}{2\pi i}\int _\gamma e^{-\zeta t}(A_\varepsilon -\zeta I)^{-1}\,d\zeta$, where $\gamma \subset\mathbb{C}$ is a contour enclosing the spectrum of $A_\varepsilon$ in the positive direction. But in \cite{MSu} only problems in a bounded domain $\mathcal{O}\subset\mathbb{R}^d$ were studied and operators under consideration were positive definite. In the case of the Dirichlet boundary condition, it was obtained that
\begin{equation}
\label{Intr exp D}
\Vert e^{-tA_{D,\varepsilon} }-e^{-tA_D^0}\Vert _{L_2(\mathcal{O})\rightarrow L_2(\mathcal{O})}\leqslant C\varepsilon (t+\varepsilon ^2)^{-1/2}e^{-ct}.
\end{equation}
The unique conceptual difference between \eqref{Intr exp} and \eqref{Intr exp D} is the behaviour at $t\rightarrow \infty$. While \eqref{Intr exp D} contains exponentially decaying factor $e^{-ct}$, we can not speak about stabilization for $A_{D,\varepsilon}$. Indeed, the difference of the operator exponentials satisfies the rough estimate $$\Vert e^{-tA_{D,\varepsilon} }-e^{-tA_D^0}\Vert\leqslant 2e^{-c_*t},$$ where $c_*>0$ is a common lower bound for $A_{D,\varepsilon}$ and $A_D^0$. In \eqref{Intr exp D}, the constant $c$ is such that $0<c<c_*$ and the constant $C$ depends on our choice of $c$ and grows as $c\rightarrow c_*$. This is caused by the used Cauchy integral representation and the behaviour of the error estimate in approximation for the resolvent $(A_{D,\varepsilon} -\zeta I)^{-1}$ for small fixed $\vert \zeta\vert$. According to the results of \cite{Su15}, the error estimate for the resolvent $(A_{\varepsilon} -\zeta I)^{-1}$ has different behaviour with respect to $\zeta$ compared to the known one for $(A_{D,\varepsilon} -\zeta I)^{-1}$. 

The \textit{goal} of the present note is to show how parabolic results from \cite{Su04,ZhPAs_parabol,Su_MMNP} can be derived from approximations for $(A_{\varepsilon} -\zeta I)^{-1}$  in $(L_2\rightarrow L_2)$- and $(L_2\rightarrow H^1)$-norms from \cite{Su15}. 
The difference between methods of the present paper and  \cite{MSu} consists of choosing the contour $\gamma$ depending on time $t$. This idea is inspired by the proof of \cite[Lemma~1]{INZ}.

\subsection*{Acknowledgement} 
The author is happy to thank the Hausdorff Research Institute for Mathematics in Bonn, Germany,  
for financial support, hospitality, and excellent working conditions during the Junior Trimester Program ,,Randomness, PDEs and Nonlinear Fluctuations.'' The author has to thank prof. T. A. Suslina for her attention to the work.

\section{Preliminaries. Known results}
\label{Section Preliminaries}

Let $\Gamma \subset \mathbb{R}^d$ be a lattice, and let $\Omega$ be the cell of the lattice $\Gamma$. By $H^1_{\mathrm{per}}(\Omega)$ we denote the subspace of matrix-valued functions from $H^1(\Omega)$ whose $\Gamma$-periodic extension belongs to $H^1_{\mathrm{loc}}(\mathbb{R}^d)$. 
For any $\Gamma$-periodic matrix-valued function $f$ we use the notation $f^\varepsilon (\mathbf{x}):=f(\varepsilon ^{-1}\mathbf{x})$, $\varepsilon >0$.  By $[f^\varepsilon]$ we denote the operator of multiplication by the matrix-valued function $f^\varepsilon (\mathbf{x})$.

In $L_2(\mathbb{R};\mathbb{C}^n)$, we consider a~matrix elliptic second order differential operator $A_\varepsilon$, $\varepsilon >0$, formally given by the expression $A_\varepsilon =b(\mathbf{D})^*g^\varepsilon (\mathbf{x})b(\mathbf{D})$. Here $g$ is a $\Gamma $-periodic $(m\times m)$-matrix-valued function, $g(\mathbf{x})>0$, $g,g^{-1}\in L_\infty$, and $b(\mathbf{D})=\sum _{l=1}^d b_lD_l$ is a first order differential operator whose coefficients $b_l$, $l=1,\dots,d$, are constant $(m\times n)$-matrices.  
The entries of the matrices $g(\mathbf{x})$ and $b_l$, $l=1,\dots,d$, are in general complex. Suppose that $m\geqslant n$ and that 
the symbol $b(\boldsymbol{\xi})=\sum _{l=1}^d b_l\xi_l$ satisfies the full rank condition: 
$\mathrm{rank}\,b(\boldsymbol{\xi})=n$, $0\neq \boldsymbol{\xi}\in\mathbb{R}^d.$ 
Or, equivalently, there exist constants $\alpha _0$ and $\alpha _1$ such that
\begin{equation*}
\alpha _0\mathbf{1}_n\leqslant b(\boldsymbol{\theta})^*b(\boldsymbol{\theta})\leqslant \alpha _1\mathbf{1}_n,\quad \boldsymbol{\theta}\in\mathbb{S}^{d-1},\quad 0<\alpha_0\leqslant \alpha _1<\infty .
\end{equation*}
Under the above assumptions, the operator $A_\varepsilon$ is self-adjoint, non-negative and strongly elliptic. The precise definition of $A_\varepsilon$ is given via the corresponding quadratic form on $H^1(\mathbb{R}^d;\mathbb{C}^n)$. 

The simplest example of the operator under consideration it the acoustics operator $A_\varepsilon =-\mathrm{div}\,g^\varepsilon(\mathbf{x})\nabla$. The operator of elasticity theory also can be written as $b(\mathbf{D})^*g^\varepsilon (\mathbf{x}) b(\mathbf{D})$, see details in \cite[Chapter 5]{BSu}.

The coefficients of the operator $A_\varepsilon$ oscillate rapidly as $\varepsilon \rightarrow 0$. The limit behaviour of its resolvent or the semigroup $e^{-tA_\varepsilon}$ is given by the corresponding function of the so-called effective operator $A^0=b(\mathbf{D})^*g^0b(\mathbf{D})$ with the constant matrix $g^0$. The definition of $g^0$ is given in terms of the $\Gamma$-periodic $(n\times m)$-matrix-valued function $\Lambda$: 
\begin{equation*}
g^0=\vert \Omega \vert ^{-1}\int _{\Omega} g(\mathbf{x})(b(\mathbf{D})\Lambda (\mathbf{x})+\mathbf{1}_m)\,d\mathbf{x},
\end{equation*}
where $\Lambda\in H^1_{\textnormal{per}}(\Omega)$ is the weak solution of the cell problem
\begin{equation*}
b(\mathbf{D})^*g(\mathbf{x})(b(\mathbf{D})\Lambda (\mathbf{x})+\mathbf{1}_m)=0,\quad \int _{\Omega }\Lambda (\mathbf{x})\,d\mathbf{x}=0.
\end{equation*}

By $S_\varepsilon$ we denote the Steklov smoothing operator acting in $L_2(\mathbb{R}^d;\mathbb{C}^m)$ by the rule
\begin{equation*}
(S_\varepsilon \mathbf{u})(\mathbf{x})=\vert \Omega\vert ^{-1}\int _\Omega \mathbf{u}(\mathbf{x}-\varepsilon \mathbf{z})\,d\mathbf{z}.
\end{equation*}
According to \cite[Lemma~1.2]{ZhPas2}, for any $\Gamma$-periodic function $f$ in $\mathbb{R}^d$ such that $f\in L_2(\Omega)$, the operator $[f^\varepsilon] S_\varepsilon$ is continuous in $L_2(\mathbb{R}^d)$, and $\Vert [f^\varepsilon] S_\varepsilon \Vert _{L_2(\mathbb{R}^d)\rightarrow L_2(\mathbb{R}^d)}\leqslant \vert \Omega\vert ^{-1/2}\Vert f\Vert _{L_2(\Omega)}$. Using this fact and the inclusion $\Lambda \in H^1_{\mathrm{per}}(\Omega)$, one can show that the so-called corrector 
\begin{equation}
\label{ell corrector}
K(\varepsilon ;\zeta):=[\Lambda ^\varepsilon] S_\varepsilon b(\mathbf{D})(A^0-\zeta I)^{-1}
\end{equation}
acts continuously from $L_2(\mathbb{R}^d;\mathbb{C}^n)$ to $H^1(\mathbb{R}^d;\mathbb{C}^n)$, and $\Vert K(\varepsilon ;\zeta)\Vert _{L_2\rightarrow H^1}=O(\varepsilon ^{-1})$ for fixed $\zeta\in\mathbb{C}\setminus\mathbb{R}_+$. The $(L_2\rightarrow H^1)$-continuity of the operator \eqref{corr time} below can be checked with the help of the same arguments.

The following result was obtained in \cite[Theorems 2.2 and 2.4]{Su15}.
\begin{theorem}[\cite{Su15}]
Let the above assumptions be satisfied. Let $\zeta\in\mathbb{C}\setminus\mathbb{R}_+$, $\phi=\mathrm{arg}\,\zeta$. Denote 
\begin{equation}
\label{c(phi)}
c(\phi)=\begin{cases}
\vert \sin\phi\vert ^{-1},\quad\phi\in (0,\pi/2)\cup (3\pi/2,2\pi),
\\
1,\quad \phi\in [\pi/2,3\pi/2].
\end{cases}
\end{equation}
Then for $\varepsilon >0$ we have
\begin{equation}
\label{Th resolvent L2}
\Vert (A_\varepsilon -\zeta I)^{-1}-(A^0-\zeta I)^{-1}\Vert _{L_2(\mathbb{R}^d)\rightarrow L_2(\mathbb{R}^d)}\leqslant C_1 c(\phi)^2\vert \zeta\vert ^{-1/2}\varepsilon .
\end{equation}
Let $K(\varepsilon ;\zeta)$ be the corrector \eqref{ell corrector}. Then for $\varepsilon >0$ we have
\begin{align}
\label{Th D res corr}
\Vert &\mathbf{D}((A_\varepsilon -\zeta I)^{-1}-(A^0-\zeta I)^{-1}-\varepsilon K(\varepsilon ;\zeta) )\Vert _{L_2(\mathbb{R}^d)\rightarrow L_2(\mathbb{R}^d)}\leqslant C_2 c(\phi)^2\varepsilon ,
\\
\label{Th res corr}
\Vert &(A_\varepsilon -\zeta I)^{-1}-(A^0-\zeta I)^{-1}-\varepsilon K(\varepsilon ;\zeta)\Vert _{L_2(\mathbb{R}^d)\rightarrow L_2(\mathbb{R}^d)}\leqslant C_3 c(\phi)^2\vert \zeta\vert ^{-1/2} \varepsilon . 
\end{align}
The constant $C_1$ depends only on $\alpha _0$, $\alpha _1$, $\Vert g\Vert _{L_\infty}$, $\Vert g^{-1}\Vert _{L_\infty}$, and parameters of the lattice $\Gamma$. The constants $C_2$ and $C_3$ depend on the same parameters and also on $m$ and $d$.
\end{theorem}

The aim of the present paper is to give another proof of the following theorem.

\begin{theorem}[\cite{Su04},\cite{ZhPAs_parabol},\cite{Su_MMNP}]
Under the above assumptions, for $\varepsilon >0$ and $t\geqslant 0$ we have
\begin{equation}
\label{Th exp}
\Vert e^{-tA_\varepsilon}-e^{-tA^0}\Vert _{L_2(\mathbb{R}^d)\rightarrow L_2(\mathbb{R}^d)}\leqslant C_4\varepsilon (t+\varepsilon ^2)^{-1/2}.
\end{equation}
Denote 
\begin{equation}
\label{corr time}
\mathcal{K}(\varepsilon ;t):=[\Lambda ^\varepsilon] S_\varepsilon b(\mathbf{D})e^{-tA^0}. 
\end{equation}
Then for $\varepsilon >0$ and $t>0$ we have 
\begin{align}
\label{Th exp corr D}
&\bigl\Vert \mathbf{D}\bigl(e^{-tA_\varepsilon}-e^{-tA^0}-\varepsilon\mathcal{K}(\varepsilon;t)\bigr)\bigr\Vert _{L_2(\mathbb{R}^d)\rightarrow L_2(\mathbb{R}^d)}\leqslant C_5\varepsilon t^{-1},
\\
\label{Th exp corr L2}
&\Vert e^{-tA_\varepsilon}-e^{-tA^0}-\varepsilon\mathcal{K}(\varepsilon;t)\Vert _{L_2(\mathbb{R}^d)\rightarrow L_2(\mathbb{R}^d)}\leqslant C_6\varepsilon t^{-1/2}.
\end{align}
The constant $C_4$ depends only on $\alpha _0$, $\alpha _1$, $\Vert g\Vert _{L_\infty}$, $\Vert g^{-1}\Vert _{L_\infty}$, and parameters of the lattice $\Gamma$. The constants $C_5$ and $C_6$ depend on the same parameters and also on $m$ and $d$.
\end{theorem}

\begin{remark}
\textnormal{
Estimate \eqref{Th exp} was announced in \textnormal{\cite[Theorem 1]{Su04}} and proved in \textnormal{\cite[Theorem 7.1]{Su07}} and, for the acoustics operator, in \textnormal{\cite[Theorem 1.1]{ZhPAs_parabol}}. For the scalar elliptic operator $A_\varepsilon =-\mathrm{div}\,g^\varepsilon (\mathbf{x}) \nabla$, where $g(\mathbf{x})$ is a~symmetric matrix with real entries, one has $\Lambda \in L_\infty$ and it is possible to replace the smoothing operator $S_\varepsilon$ in the corrector by the identity operator. In this case, estimate \eqref{Th exp corr D} was obtained in \textnormal{\cite[Theorem 1.3]{ZhPAs_parabol}}. For the matrix elliptic operator,  
$(L_2\rightarrow H^1)$-approximation for $e^{-tA_\varepsilon}$  was proved in \textnormal{\cite[Theorem 11.1]{Su_MMNP}} (but with another smoothing operator in the corrector). 
}
\end{remark}

\section{New proof}

Using the Riesz–Dunford functional calculus, we  represent the operator exponential $e^{-tA_\varepsilon}$ as an integral: 
\begin{equation*}
e^{-tA_\varepsilon}=-\frac{1}{2\pi i}\int _\gamma
e^{-\zeta t}(A_\varepsilon -\zeta I)^{-1}\,d\zeta .
\end{equation*}
Here $\gamma$ is a suitable contour in the complex plain enclosing the spectrum $\sigma (A_\varepsilon)\subset [0,\infty)$ in the positive direction. 
One can choose $\gamma =\widehat{\gamma }\cup\widetilde{\gamma}$ with
\begin{align*}
\widehat{\gamma}&=\lbrace \zeta \in \mathbb{C} : \zeta =e^{i\phi}, \pi/4\leqslant \phi \leqslant  7\pi/4\rbrace ,
\\
\widetilde{\gamma}&=\lbrace \zeta \in\mathbb{C} : \zeta =re^{i\pi/4}, r\geqslant 1\rbrace
\cup \lbrace \zeta \in\mathbb{C} : \zeta =re^{i 7\pi/4}, r\geqslant 1\rbrace .
\end{align*}
But we take the contour depending on $t>0$, shrinking this contour $\gamma$ in $t$ times: $\gamma _t=t^{-1}\gamma=\lbrace \zeta \in \mathbb{C} : \zeta= t^{-1}\eta, \eta \in\gamma\rbrace$. 
Applying these arguments to the operators $e^{-tA_\varepsilon}$ and $e^{-tA^0}$ and changing variable, we get
\begin{equation}
\label{exp identity}
\begin{split}
e^{-tA_\varepsilon}-e^{-tA^0}&=-\frac{1}{2\pi i}\int _{\gamma _t} e^{-\zeta t}\left( (A_\varepsilon -\zeta I)^{-1}-(A^0-\zeta I)^{-1}\right)\,d\zeta 
\\
&= -\frac{1}{2\pi i t}\int _\gamma e^{-\eta}\left( (A_\varepsilon -t^{-1}\eta I)^{-1}-(A^0-t^{-1}\eta I)^{-1}\right)\,d\eta .
\end{split}
\end{equation}
Recall notation \eqref{c(phi)}. 
Using \eqref{Th resolvent L2} and taking into account that $c(\phi _t)\leqslant 2^{1/2}$, where $\eta\in\gamma$ and $\phi_t:=\mathrm{arg}\,(t^{-1}\eta)$, for $t>0$ we have
\begin{equation}
\label{exp-exp preliminary L2}
\begin{split}
\Vert  e^{-tA_\varepsilon}-e^{-tA^0}\Vert _{L_2(\mathbb{R}^d)\rightarrow L_2(\mathbb{R}^d)}
\leqslant \frac{C_1\varepsilon}{\pi t}\int _\gamma \vert e^{-\eta}\vert \vert t^{-1}\eta\vert ^{-1/2}\vert\,d\eta\vert
\leqslant 
\frac{C_1\varepsilon}{\pi t^{1/2}}\int _\gamma \vert e^{-\eta}\vert \vert\,d\eta\vert.
\end{split}
\end{equation}
The integral here is understood as a contour integral of the first kind. Let us estimate it:
\begin{equation}
\label{c appears}
\frac{1}{\pi}\int _\gamma \vert e^{-\eta}\vert \vert\,d\eta\vert
\leqslant \frac{1}{\pi}\int _{\pi/4}^{7\pi /4}e^{-\cos \phi}\,d\phi + \frac{2}{\pi}\int _1^\infty e^{-r/\sqrt{2}}\,dr\leqslant 3e/2 + 2^{3/2}\pi ^{-1}e^{-1/\sqrt{2}}=:\mathfrak{c}.
\end{equation}
Thus, $\Vert  e^{-tA_\varepsilon}-e^{-tA^0}\Vert _{L_2(\mathbb{R}^d)\rightarrow L_2(\mathbb{R}^d)}
\leqslant \mathfrak{c}C_1\varepsilon t^{-1/2}$. 
Obviously, for $t\geqslant 0$ the left-hand side of \eqref{exp-exp preliminary L2} does not exceed $2$. Since $\min\lbrace 2; \mathfrak{c}C_1\varepsilon t^{-1/2}\rbrace \leqslant C_4 \varepsilon (t+\varepsilon ^2)^{-1/2}$, where $C_4=2^{1/2}\max\lbrace 2 ;\mathfrak{c}C_1\rbrace$, we arrive at estimate \eqref{Th exp}.

To prove $(L_2\rightarrow H^1)$-approximation in the same manner, we need an identity for the correctors \eqref{ell corrector} and \eqref{corr time}. So, we act by the operator 
 $[\Lambda^\varepsilon]S_\varepsilon b(\mathbf{D})$ from the left to the both sides of the contour integral representation for the exponential $e^{-tA^0}$. Since the operator $[\Lambda^\varepsilon]S_\varepsilon b(\mathbf{D})$ is closed, we can move it across  the integral sign. Thus, 
\begin{equation}
\label{identity exp corrector}
\mathcal{K}(\varepsilon;t)=-\frac{1}{2\pi i}\int _{\gamma_t} e^{-\zeta t}K(\varepsilon ;\zeta)\,d\zeta .
\end{equation}
Similarly to the proof of estimate \eqref{Th exp}, relations \eqref{Th D res corr}, \eqref{exp identity}, \eqref{c appears}, and \eqref{identity exp corrector} imply estimate \eqref{Th exp corr D} with the constant $C_5:=\mathfrak{c}C_2$. Estimate \eqref{Th exp corr L2} follows from \eqref{Th res corr} on the same way, $C_6:=\mathfrak{c}C_3$.

\section{Discussion}

Since we derive the parabolic estimates from the elliptic ones, the achievement of the present paper can be interpreted  as a quantitative Trotter-Kato like result in homogenization context. For derivation of hyperbolic results from elliptic ones, see preprint \cite{M}.

The author believes that the used technique may be useful for positive definite operators after refinement of the known resolvent estimates near the lower edge of the spectrum.

\end{document}